\documentclass[10pt,reqno]{amsart}

\usepackage{amssymb,amsmath}

\def\bl{\begin{lemma}}
\def\el{\end{lemma}}
\def\bth{\begin{theorem}}
\def\eth{\end{theorem}}
\def\bc{\begin{corollary}}
\def\ec{\end{corollary}}
\def\bcj{\begin{conjecture}}
\def\ecj{\end{conjecture}}
\def\bpr{\begin{proposition}}
\def\epr{\end{proposition}}
\def\bde{\begin{definition}}
\def\ede{\end{definition}}
\def\E{\mathbb{E}}
\def\N{\mathbb{N}}
\def\Pr{\mathbb{P}}
\def\QED{\hfill\qedsymbol}

\newcommand{\br}{\mbox{\rm br}}

\newcommand{\pr}{probability }
\newcommand{\crit}{critical }
\newcommand{\anch}{\iota^*}

\newcommand{\jour}[5]{#1. #2. {\it #3} {\bf #4} #5}
\newcommand{\conf}[4]{#1. #2. {\it #3}, #4}
\newcommand{\book}[4]{#1. {\it #2}. #3, #4.}

\newcommand{\toap}[3]{#1. #2. {\it #3}, to appear.}
\newcommand{\subm}[2]{#1. #2. Submitted.}

\newcommand{\dist}{\mbox{\rm dist}}
\newcommand{\be}{\begin{eqnarray}}
\newcommand{\ee}{\end{eqnarray}}

\newcommand{\Z}{{\mathbb Z}}

\newcommand{\T}{{\mathcal T}}
\newcommand{\F}{{\mathcal F}}

\newcommand{\RR}{{\mathcal R}}

\newtheorem{theorem}{Theorem}[section]
\newtheorem{definition}{Definition}[section]
\newtheorem{lemma}[theorem]{Lemma}

\newtheorem{proposition}[theorem]{Proposition}

\theoremstyle{definition}
\numberwithin{equation}{section}

\input epsf.sty

\begin{document}

\title{Bootstrap Percolation on Infinite Trees
and non-amenable groups}

\author{J\'ozsef Balogh}
\address{Department of Mathematics, The Ohio State University, 231 W
18th Ave, Columbus, OH 43235}
\email{jobal@sol.cc.u-szeged.hu\qquad www.math.ohio-state.edu/\~{}jobal}

\author{Yuval Peres}
\address{Departments of Statistics and Mathematics, 367 Evans Hall,
University of California, Berkeley, CA 94720}
\email{peres@stat.berkeley.edu\qquad
www.stat.berkeley.edu/\~{}peres}
\author{G\'abor Pete}
\address{Department of Statistics, 367 Evans Hall, University of
 California, Berkeley, CA 94720}
\email{gabor@stat.berkeley.edu\qquad www.stat.berkeley.edu/\~{}gabor}

\thanks{Our work was partially supported by NSF grants DMS-0302804 (Balogh),
  DMS-0104073 and DMS-0244479 (Peres, Pete), and OTKA (Hungarian
  National Foundation for Scientific Research) grants T34475 (Balogh) and
  T30074 (Pete).} 

\date{April 19, 2005. \smallskip}

\begin{abstract}
Bootstrap percolation on an arbitrary graph has a random initial
configuration, where each vertex is occupied with probability $p$,
independently of each other, and a deterministic spreading rule
with a fixed parameter $k$: if a vacant site has at least $k$
occupied neighbors at a certain time step, then it becomes
occupied in the next step. This process is well-studied on $\Z^d$;
here we investigate it on regular and general infinite trees and on 
non-amenable Cayley graphs. The \crit \pr is the infimum of those 
values of $p$ for which the process achieves complete occupation with 
positive probability. On trees we find the following
discontinuity: if the branching number of a tree is strictly
smaller than $k$, then the \crit \pr is 1, while it is $1-1/k$ on
the $k$-ary tree. A related result is that in any rooted tree $T$
there is a way of erasing $k$ children
of the root, together with all their descendants, and repeating
this for all remaining children, and so on, such that the
remaining tree $T'$ has branching number $\br(T')\leq
\max\{\br(T)-k,\,0\}$. We also prove that on any $2k$-regular
non-amenable graph, the \crit \pr for the $k$-rule is strictly
positive.
\end{abstract}

\maketitle

\section{Introduction and results}\label{intro}

Consider a countable, connected, locally finite graph $G=G(V,E)$, with two 
possible states for each site in the vertex set $V$: vacant (0) or occupied
(1). Start with a
configuration picked according to the product Bernoulli measure $\Pr_p$,
i.e.~each site is occupied randomly and
independently with \pr $p$. Then fix a parameter $k$, and consider
the following deterministic spreading rule: if a vacant site has at least $k$
occupied neighbors at a certain time step, then it becomes occupied
in the next step. This process is called {\bf bootstrap percolation}.
{\bf Complete occupation} is the event that every vertex becomes
occupied during the process. The main problem is to determine the
{\bf \crit probability} $p(G,k)$ for complete occupation: for infinite graphs
$G$ this is the infimum of the initial probabilities $p$ that make
$\Pr_p$(complete occupation$)>0$. This model has a rich history in
statistical physics, mostly on $G=\Z^d$ and finite boxes; we will give
some references later. 

For infinite trees the most important characteristic of growth is
the {\bf branching number} $\br(T)$ of the tree, see \cite{ly:perc90}
or \cite{lpbook}. It is defined as the supremum of real numbers
$\lambda\geq 1$ such that $T$ admits a positive flow from the root to
infinity, where on every edge $e\in E(T)$, the flow is bounded by
$\lambda^{-|e|}$, and $|e|$ denotes the number of edges (including
$e$) on the path from $e$ to the root. This supremum does not depend
on the root, and remains unchanged if we modify a finite portion of
the tree. Two basic examples are $\br(T_k)=k$ for the
$(k+1)$-regular tree, and $\br(T_\xi)=\E\xi$ a.s.~given non-extinction
for the Galton-Watson tree $T_\xi$ with offspring distribution
$\xi$. For finite trees, the branching number is 0.

On $T_k$, $k$-neighbor bootstrap percolation has
$p(T_k,k)=1-1/k$, see (\ref{far}) in Proposition \ref{regtree} below.
In contrast, we have the following:

\bth\label{maintree}
Let $T$ be an infinite tree. If $\br(T)<k$, then $p(T,k)=1$.
\eth

The above results show a somewhat surprising discontinuity of the function
\be\label{function}
f_k(b):=\inf\{p(T,k)\,:\,\br(T)\leq b,\ T\text{ has bounded degree}\}
\ee
at the value $b=k$. If we omit the condition of bounded degree, the
discontinuity is even sharper: it is easy to construct a tree with
$\br(T)=k$ and $p(T,k)=0$. A possible explanation of this
discontinuity is given by Theorem \ref{pruning} below.

For regular trees we give an equation for the \crit probability,
from which the actual value is more-or-less computable.

\bpr\label{regtree}
Let $2\leq k\leq d$. The \crit \pr $p(T_d,k)$ is the supremum of all $p$ for
which the equation 
\be\label{binom}
\Pr\bigl({\rm Binom}(d,(1-x)(1-p))\leq d-k\bigr)=x
\ee
has a real root $x\in(0,1)$. In particular, for any constant
$\gamma\in [0,1]$ and a sequence of integers $k_d$ with
$\lim_{d\to\infty} {k_d}/{d}=\gamma$,
\be\label{middle}
\lim_{d\to\infty}p(T_d,k_d)=\gamma.
\ee
Furthermore, for the extreme values of the parameter $k$,
\be\label{far}
p(T_d,d)=1-\frac{1}{d}\mbox{\qquad{\rm and}\qquad}
p(T_d,2)=1-\frac{(d-1)^{2d-3} }{ d^{d-1}(d-2)^{d-2}}\sim \frac{1}{2d^2}.
\ee
\epr

\medskip

There is a generalization of a weaker form of Theorem \ref{maintree}.
For this we first have to introduce the following simple notion, which
will also be central to our proofs.

\bde\label{fort} A finite or infinite connected subset $F\subseteq V$
of vertices is called a {\bf $k$-fort} if each $v\in F$ has outdegree
$\deg_{V\setminus F}(v)\leq k$. Here $\deg_H(v)=|\{w\in H:(v,w)\in E\}|$,
for any $H\subseteq V$.
\ede

A key observation is that the failure of complete occupation by the
$k$-neighbor rule is equivalent to the existence of a vacant
$(k-1)$-fort in the initial configuration.

\bth\label{pruning} Let $T$ be an infinite tree. Then every vertex
$x\in T$ is contained in a $k$-fort $F$ with $\br(F)\leq
\max\{\br(T)-k,\,0\}$.
\eth

This means that after fixing any vertex as the root, we can erase
$k$ children of it, together with all their descendants, and can
repeat this for all the remaining children, and so on, so that
this pruning process results in a required subtree $F$. It is
interesting to note that the natural idea of pruning off the $k$
subtrees with the largest branching numbers at each generation
does not work in general.

For $\br(T)<k$ we get a $(k-1)$-fort with $\br(F)<1$, which can
happen only if $F$ is finite, so $\br(F)=0$. In fact, in Theorem
\ref{maintree} we prove that there are infinitely many finite
$(k-1)$-forts of bounded size, which implies $p(T,k)=1$. The
impossibility of $0<\br(F)<1$ might be viewed as the reason for
the discontinuity of $f_k(b)$ at $b=k$, though we do not actually
know continuity at other points. See Section \ref{open} for more
discussion and open problems.
\medskip

An infinite graph $G$ has the {\bf anchored expansion property} if
for some fixed vertex $o\in V(G)$, the {\bf anchored Cheeger constant}
is positive:
\be\label{anchchee}
0<\anch(G):=\liminf\left\{\frac{|\partial_e S|}{|S|}:o\in S\subset V,
\ S \mbox{ is finite and connected}\right\},
\ee
where $\partial_e S$ is the set of edges in $E(G)$ with exactly one
endpoint in $S$. It is easy to see that the value of $\anch(G)$ does
not depend on the vertex $o$. This notion is implicit in \cite{tho},
and was defined explicitly by \cite{bls:pertu}. For transitive graphs
(such as Cayley graphs of finitely generated infinite groups) it
coincides with the more
familiar but less robust concept of non-amenability, where the infimum
is taken over all finite connected subsets $S$. For background on
non-amenability see \cite{lpbook} or \cite{ly:phase}, and on anchored
expansion \cite{hss} or \cite{vir:anch}.

\bth\label{anchthm}
Let $G_d$ be a $d$-regular graph. If $\anch(G_d)+2k>d$, then
$p(G_d,k)>0$. In particular, if $G_d$ has the anchored expansion
property, then $p(G_d,\lceil d/2\rceil)>0$.
\eth

This result is sharp in the sense that there exists a 6-regular
non-amenable Cayley graph $G_6$ with $p(G_6,2)=0$, see Section
\ref{anchored}.
We will pose a possible characterization of amenability
in Section \ref{open}.
\medskip

The issue of positivity of the \crit \pr is simpler for the case of 
trees. For this, let us denote by $q(G,k)$ the 
infimum of initial probabilities for which, following the $k$-neighbor 
rule on $G$, there will be an infinite connected component of occupied 
vertices in the final configuration with positive probability. 
Clearly, $q(G,k)\leq p(G,k)$. 

\bpr\label{positive}
For any integer $d$, and $k\geq 2$, if $T$ is an infinite tree with 
maximum degree $d+1$, then $p(T,k)\geq q(T_d,k)>0$.
\epr

The first inequality of this proposition follows immediately from 
viewing $T$ as a subgraph of $T_d$. The positivity of the \crit \pr 
$q(T_d,k)$ will be proved using our proof of Proposition 
\ref{regtree} and an idea from \cite{how}. 
\medskip

Bootstrap percolation was first defined in the statistical physics
literature in \cite{cha79}, where the formulae of (\ref{far}) were
given. A variant of the model appeared in \cite{cha82}. The
problem of complete occupation on $\Z^2$ was solved by
\cite{vEnt}.
Schonmann proved \cite{sch92} that the \crit \pr $p(\Z^d,k)$
for bootstrap percolation is 0 for $k\leq d$ and is 1 for $k>d$.
The process can also be considered on finite graphs, see
e.g.~\cite{aizleb}, \cite{babo:hyp} and \cite{hol:meta}. A short
recent physics survey is \cite{adlev}. Bootstrap percolation also 
has connections to the dynamics of the Ising model at zero
temperature; see \cite{fss} for $\Z^d$, and \cite{how} for
$T_2$. 
\medskip

We conclude this introduction by some basic observations.

If a graph $G$ satisfies $\Pr_p($complete occupation with the 
$k$-rule$)\in\{0,1\}$ for all $p\in[0,1]$, and so $\Pr_p($complete 
occupation of $G)=1$ for any $p>p(G,k)$, then we will say that 
the {\bf 0-1 law holds} for $G$ with the $k$-rule.

For example, if the orbit of each vertex under the automorphism group 
of $G$ is infinite, then the product \pr measure of the initial configuration
is ergodic \cite[Proposition 6.3]{lpbook}, while complete occupation is
an invariant property, hence it has probability 0 or 1. Furthermore,
if there is a finite $(k-1)$-fort in such a $G$, we immediately have
infinitely many copies of this, so $p(G,k)=1$. On the other hand:

\bl\label{nofinite} If there are no finite $(k-1)$-forts in a graph
$G$, then  $p(G,k)\leq 1-p_c(G)$, where $p_c(G)$ denotes the \crit \pr
for standard site  percolation on $G$.  \el

\proof In the case of no complete occupation, the vacant $(k-1)$-fort
has to be infinite, thus we have an infinite connected vacant
component in the  initial configuration. To have this event with
positive probability, the density of initial vacant sites has to be at
least the \crit \pr $p_c(G)$.\QED
\smallskip

Therefore, if $p_c(G)>0$ holds for a graph without finite
$(k-1)$-forts, which is usually the case (e.g.~if the degrees of
vertices  are bounded, see \cite[Prop.~6.9]{lpbook}), then $p(G,k)<1$.
For instance, on any tree $T$ we have $p_c(T)=1/\br(T)$, as was shown
in \cite{ly:perc90}.

We will say that a graph $G$ is {\bf uniformly bigger} than a graph 
$H$ if every vertex of $G$ is contained in a subgraph of $G$ that is 
isomorphic to $H$. 

\bl\label{monotone} {\bf (Monotonicity)} If a graph $G$ is uniformly 
bigger than $H$, and $H$ satisfies the 0-1 law for some $k$-rule, 
then we have $p(G,k)\leq p(H,k)$.
\el

\proof For any $p>p(H,k)$, any fixed vertex $v$ of $G$ becomes occupied
almost surely, because of the copy of $H$ containing $v$. There are
countably many vertices of $G$, so we have $\Pr_p($complete occupation 
of $G)=1$ with this $p$.\QED
\smallskip

In particular, if $T$ is a tree with maximal degree $d+1$ and it
satisfies the 0-1 law, then we get $p(T,k)\geq p(T_d,k)$. 
Proposition \ref{positive} is a generalization of this fact. We thank 
\'Ad\'am Tim\'ar for pointing out the importance of considering $q(T_d,k)$
for the generalization.  

\section{Regular trees}\label{regulartrees}

\noindent{\it Proof of Proposition \ref{regtree}.} Consider the
$(d+1)$-regular tree $T_d$, and fix $2\leq k\leq d$. This tree has no
finite $(k-1)$-forts, and it is easy to see that any infinite fort of
it contains a complete $(d+2-k)$-regular subtree. Hence, unsuccessful
complete occupation for the $k$-rule is equivalent to the existence of
a $(d+2-k)$-regular vacant subtree in the initial configuration.

Note that complete occupation on $T_d$ obeys the 0-1 law. So
incomplete occupation has \pr $1$ if and only if a fixed origin is
contained in a $(d+2-k)$-regular vacant subtree with positive
probability. Now a simple use of Harris' inequality, see \cite[Section
6.2]{lpbook}, gives that this is equivalent to having the following
event with positive probability: a $d$-ary tree, rooted at the fixed
origin that is declared to be vacant, has a vacant $(d+1-k)$-ary
subtree starting from the same root. Therefore, we need to determine
when the connected component of vacant sites of the root, which is a
random Galton-Watson tree with offspring distribution Binom$(d,1-p)$,
contains a $(d+1-k)$-ary subtree with positive probability.  If the
\pr of {\em not} having such a subtree is denoted by $y=y(p)$, then each
of the $d$ children of the root has probability $1-p$ to be vacant, and
given this event, has probability $1-y$ to be the root of a vacant
$(d+1-k)$-ary subtree. Therefore, $y$ clearly satisfies the equation
(\ref{binom}), i.e. it is a fixed point of the function
\be
x \mapsto B_{d,k,p}(x)&:=&\Pr\bigl({\rm Binom}(d,(1-x)(1-p))\leq d-k\bigr)\nonumber\\
&=&\sum_{j=0}^{d-k}{d\choose j}(1-x-p+xp)^j(p+x-px)^{d-j}.\nonumber
\ee
One fixed point in $[0,1]$ is $x=1$; we are going to show that $y$
is actually the smallest one in $[0,1]$. It is easy to see that
$$
\frac{\partial}{\partial x}B_{d,k,p}(x)=d(1-p) \Pr\bigl({\rm
Binom}(d-1,(1-x)(1-p)) = d-k \bigr),
$$
which is positive for $x\in [0,1)$, with at most one extremal point (a
maximum) in $(0,1)$. Thus $B_{d,k,p}(x)$ is a monotone increasing
function with $B_{d,k,p}(0)>0$ and with at most one inflection
point in $(0,1)$. If $y_n$ denotes the probability that the
required vacant subtree does not even reach the $n$th level below
the root, then $y_0=0$, $y_{n+1}=B_{d,k,p}(y_n)$, and $y_n\to y$.
On the other hand, the sequence $y_n$ clearly approaches the
smallest fixed point of $B_{d,k,p}(x)$, which so coincides with
$y$. Thus,
the infimum of the probabilities $p$ for which equation
(\ref{binom}) has no positive real root $x<1$ is indeed the \crit
\pr $p(T_d,k)$.

If $\lim_{d\to\infty} {k_d}/{d}=\gamma$, then for any fixed $p$ and
$x$, by the Weak Law of Large Numbers:
$$
B_{d,k_d,p}(x)=\Pr\left(\frac{{\rm Binom}(d,(1-x)(1-p))}{d}
\leq\frac{d-k_d}{d}\right)\to\begin{cases} 1,&\mbox{if\
}(1-x)(1-p)<1-\gamma\\ 0,&\mbox{if\ }(1-x)(1-p)>1-\gamma,
\end{cases}
$$
as $d\to\infty$. Solving the equation $(1-x)(1-p)=1-\gamma$ for
$x$ gives a \crit value $x_c=(\gamma-p)/(1-p)$. Thus for
$p>\gamma$ we have $\lim_{d\to\infty}B_{d,k_d,p}(x)\to 1$ for all
$x\in [0,1]$, while for large enough $d$, $B_{d,k,p}(x)$ is 
convex in $[0,1]$, so there is no positive root $x<1$ of 
$B_{d,k_d,p}(x)=x$. On the other hand, for $p<\gamma$ there 
must be a root $x=x(d)$ for large enough $d$,
clearly satisfying $\lim_{d\to\infty} x(d)= x_c$. These prove
(\ref{middle}).

The first equality of (\ref{far}) follows immediately from
(\ref{binom}). The second equality can be deduced by a standard
calculus argument from our above formula for the first derivative
of $B_{d,k,p}(x)$.
\QED
\smallskip

\noindent{\bf Remark 1.} We will use later that the extinction
probability $y(p)$ introduced in the above proof satisfies $y(p)\to 0$
as $p\to 0$. This follows from the facts that $B_{d,k,0}(x)<x$
for $x>0$ small enough, that the functions $B_{d,k,p}(x)$ converge 
uniformly to $B_{d,k,0}(x)$ as $p\to 0$, and that $B_{d,k,p}(0)>0$.
\smallskip 

\noindent{\bf Remark 2.} A Galton-Watson tree with offspring
distribution Binom$(d,1-p)$ can contain a $(d+1-k)$-ary subtree only
if its mean is $(1-p)d\geq d+1-k$. Thus $p(T_d,k)\leq \frac{k-1}{d}$
follows immediately.
\smallskip

\noindent{\bf Remark 3.} The problem of finding regular subtrees
in certain Galton-Watson trees was first considered in
\cite{chchd}, where the formula of (\ref{far}) for $p(T_9,2)$ was
used. For general Galton-Watson processes, see \cite{pakdek}. From
($\ref{middle}$) it follows that the critical mean value for a
binomial offspring distribution to produce an $N$-ary subtree in
the Galton-Watson tree is asymptotically $N$. In \cite{pakdek} it
was shown that this critical mean value is $\sim eN$ for a
geometric offspring distribution, and $\sim N$ for a Poisson
offspring. An interesting feature of these phase transitions is
that unlike the case of usual percolation $N=1$, for $N\geq 2$ the
\pr of having the $N$-ary subtree is already positive at
criticality. For bootstrap percolation this means that the \pr of
complete occupation is still 0 at $p=p(T_d,k)$ if $k<d$.
\smallskip

\noindent{\bf Remark 4.} The \crit \pr $p(T,k)$ can be computed also 
for quasi-transitive (periodic) trees and Galton-Watson 
trees, as we will see for example in Section \ref{open}.
\medskip

\noindent{\it Proof of Proposition \ref{positive}.} To prove the
second inequality, $q(T_d,k)>0$, we will first show that for any
non-backtracking path $v_0,v_1,\dots,v_n$ in $T_d$,  
\be\label{path}
\Pr_p\bigl(\{v_0,\dots,v_n\}
\mbox{ does not intersect any vacant }
(k-1)\mbox{-fort}\bigr)\leq [1-z(p)^2]^{\lfloor n/2\rfloor}, 
\ee
where $z(p)$ is the \pr that an infinite rooted tree $Z$ with $d-1$ 
children at the root, and $d$ children everywhere else, has a $(d+1-k)$-ary 
vacant subtree containing the root. Before proving (\ref{path}), note that 
$z(p)\to 1$ as $p\to 0$. This follows easily from the fact $y(p)\to 0$ 
of Remark 1 above.

To prove (\ref{path}), for each $v_i$ consider a copy $Z_i$ of $Z$
inside $T_d$, rooted at $v_i$, disjoint from the path $v_0,\dots,v_n$. Then
the subtrees $Z_i$ are also disjoint from each other. The $(d+1-k)$-ary
vacant subtrees rooted at $v_{2i}$ and $v_{2i+1}$ inside $Z_{2i}$ and
$Z_{2i+1}$ join together to give a vacant $(d+2-k)$-regular tree, i.e.~a
$(k-1)$-fort, inside $T_d$.  The probability that this does not happen for any
of the pairs $v_{2i}, v_{2i+1}$ is exactly the RHS of (\ref{path}).

The number of different paths $v_0,v_1,\dots,v_n$ from a fixed vertex $v_0=o$ 
is $(d+1)d^{n-1}$. Therefore, if $p$ is so small that $\sqrt{1-z(p)^2}<1/d$,
then the \pr that there is at least one such path that does not intersect any vacant $(k-1)$-forts in the initial configuration is exponentially small in $n$. By the 
Borel-Cantelli lemma, any infinite non-backtracking path from $o$ eventually
intersects a vacant $(k-1)$-fort almost surely, hence the bootstrap
percolation process will not be able to form an infinite occupied cluster 
containing $o$. There are countably many possible $o$ vertices in $T_d$, 
so we have the same for all vertices with \pr 1. 
Thus $q(T_d,k)\geq p$ with the above small $p>0$.
\QED

\section{General infinite trees}\label{generaltrees}

To start our discussion of the connection between branching number
and bootstrap percolation, let us prove a simple combinatorial
lemma, which implies Theorem \ref{pruning} for the special case of
$\br(T)<k$, but is not yet enough to prove Theorem \ref{maintree}.

\bl\label{red}{\bf (Red Lemma)} If some vertex $x$ of a tree $T$ is
not contained in any finite $(k-1)$-fort, then there is a $k$-ary subtree
containing $x$.  \el

\begin{figure} \centerline{\epsfxsize=2.5 in
 \epsffile{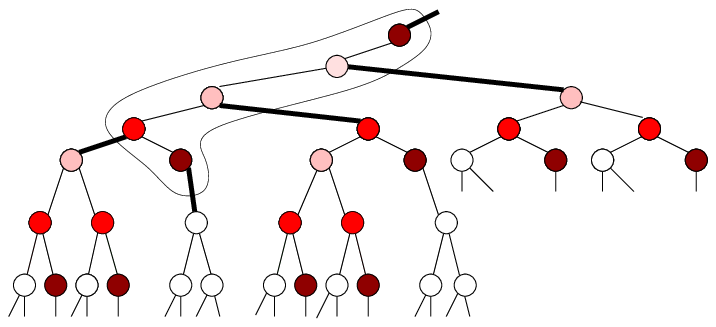} \hskip 0.5 in  \epsfxsize=1.4 in
 \epsffile{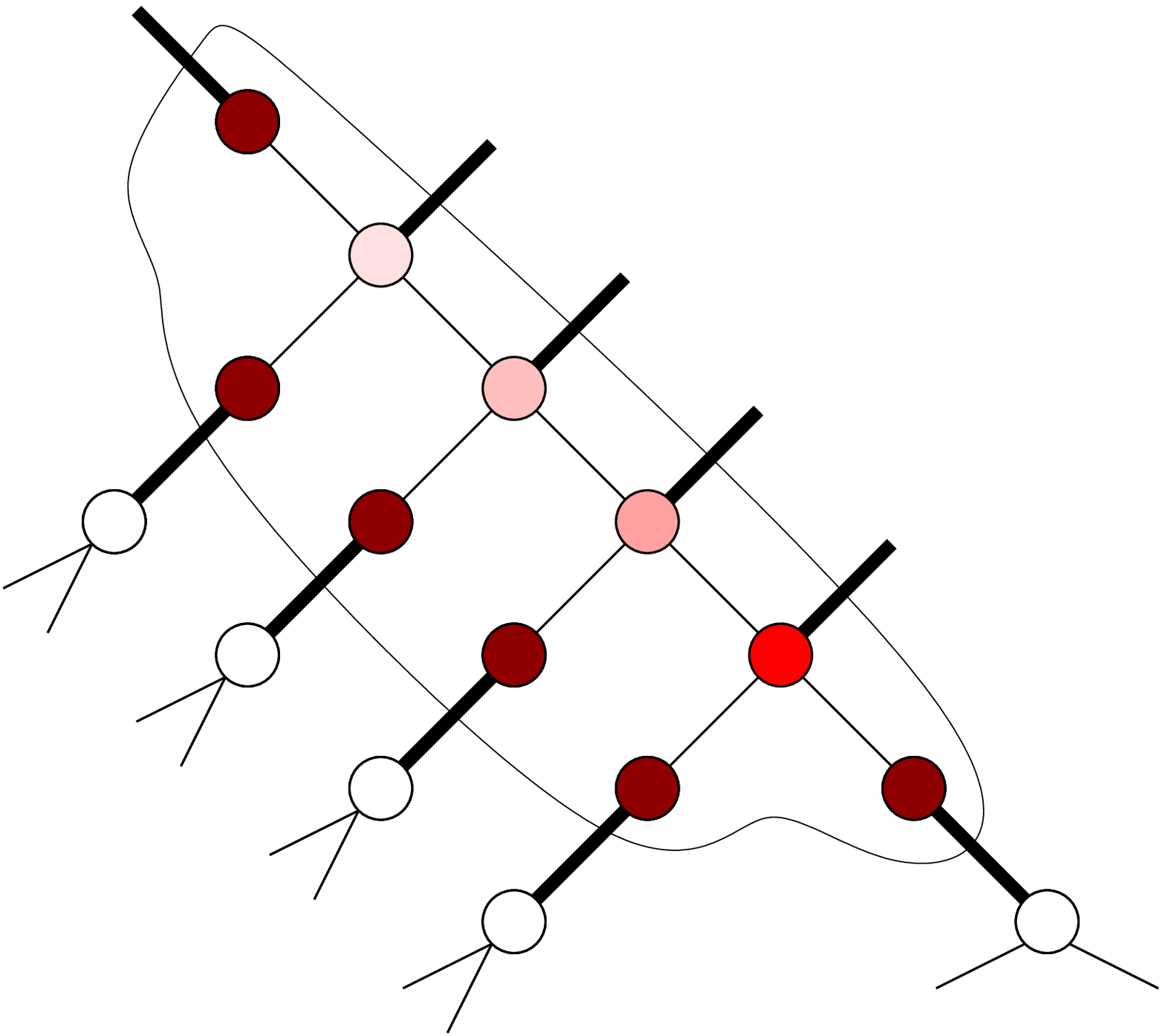}}\caption{Some finite
 1-forts.}\label{1-forts} \end{figure}

\proof Consider the tree $T$ as rooted at the vertex $x$. First color
red all vertices with at most $k-1$ children. In the second step,
color red each vertex with at most $k-1$ non-red children, and repeat
this over and over again; see Figure \ref{1-forts}. 
In the limiting final coloring, if the root $x$ is red,
then it obtained its color in a finite number of steps, so  there is a
finite set $F$ of vertices such that $x$ becomes red even if we fix
all the vertices outside $F$ to be uncolored forever. If we take this
$F$ to be minimal, then it is a finite connected subtree of $T$, with
all leaves painted red in the first step, and all vertices becoming
eventually red. But now, this red $F$ is clearly a finite $(k-1)$-fort
in $T$, contradicting the choice of $x$. Therefore, $x$ is not red in
the final limiting coloring. This means it has at least $k$ non-red
children, and each of these children also has at least $k$ non-red
children, and so on. Hence, the non-red component of $x$ in $T$
contains the $k$-ary subtree we wanted.\QED
\smallskip

Theorem \ref{maintree} is not obvious from this lemma because if
we do not forbid all finite $(k-1)$-forts, but only the appearance
of too many small ones, we can already get $p(T,k)<1$, while we
can have vertices with at most $k-1$ children lying close to each
other (see the tree on the right in Figure \ref{1-forts}), so
$\br(T)<k$ could possibly occur, as well. Thus we need a
quantitative version of Lemma \ref{red}. For this, fix a root $o$
for $T$, and denote by $L_r(x)$ the set of the vertices $y$ from
which the shortest path to $o$ contains $x$, and
$\dist(o,y)=\dist(o,x)+r$.

\bl\label{blue} {\bf (Blue Lemma)} Let $x$ be a vertex with
 $|L_R(x)|<(k-1)k^{R-1}$ for  some positive integer $R$. Then
 $\cup_{0\le r \le R} L_r(x)$ contains  a $(k-1)$-fort of $T$.
\el

\proof Let $x$ be a vertex of the tree satisfying the conditions
of the lemma, and label its level by $0$. Color a vertex on level
$R$ blue, if it has at least $k$ children. In general, color a
vertex on level $r$ blue, if it has at least $k$ blue children
($r<R$).  The vertex $x$ is definitely not blue, otherwise
$|L_R(x)|\geq k^R$ would hold. Moreover, $x$ has at most $k-2$
blue children, since $|L_R(x)|< (k-1)k^{R-1}$.  If $x$ has less
than $k-1$ children, then it is a fort by itself. Otherwise, the
non-blue component containing $x$ contains at least $2$ vertices.
We claim that the non-blue connected component containing $x$ is a
$(k-1)$-fort. First of all, $x$ has outdegree at most $k-1$,
counting  its mother and its possible $k-2$ blue children. Any
other vertex from this set has a non-blue mother, and being
non-blue means that it has at most $k-1$ blue neighbors.  A
non-blue vertex in this component in the level $R$ has at most
$k-1$ children, and its mother is not blue.\QED
\smallskip

 Note that almost the same argument for $x=o$ gives that already
 $|L_R(o)|<k^R$ implies a $(k-1)$-fort inside $\cup_{0\leq r\leq
 R}L_r(o)$.  (The reason for the strengthening is simple:
 $o$ does not have a mother.)
\smallskip

\noindent{\it Proof of Theorem \ref{maintree}.} We will prove that
if there is no $(k-1)$-fort with at most $N$ vertices,  then
$\br(T)\geq k-\frac{2k\log k}{\log N}$. This suffices because destroying
a finite number of forts of size at most $N$ does not affect
$\br(T)$, so $\br(T)<k$ will imply the existence of infinitely
many $(k-1)$-forts of bounded size, which shows $p(T,k)=1$.

Any leaf of $T$ would be a $(k-1)$-fort with one vertex, so there are no
leaves, and we have $|L_r(x)|\leq |L_{r+1}(x)|$ for any $r$ and $x$. Hence,
 the fort that the Blue Lemma finds for us has less than $Rk^R<k^{2R}$
 vertices. Thus having no $(k-1)$-forts of size at most $N$ implies
 that $|L_r(x)|\geq (k-1)k^{r-1}$ and $|L_r(o)|\geq k^r$ for every
 $r\leq R=\frac{\log N}{2\log k}$. We will prove that if $\lambda>0$ is
 such that
\be\label{capacity}
\lambda^R\leq (k-1)k^{R-1},
\ee
 then $\br(T)\geq\lambda$. For example, if $\lambda=k-\frac{c_k}{\log
 N}$, then $\lambda^R\leq N^{1/2}\exp(-\frac{c_k}{2k\log k})$, while
 $(k-1)k^{R-1}=\frac{k-1}{k}N^{1/2}$, so $c_k=2k\log k$ is good for
 any $k\geq 2$.

 Now we have to show that if the capacity of an edge $e$ is
 $\lambda^{-|e|}$, then the network admits a positive flow from the
 root to infinity. Start the flow with an amount $k^{-R}$ at the root
 $o$.  On level $R$ there are at least $k^R$ vertices; divide the
 initial amount equally among them, and build the flow from $o$ to
 $L_R(o)$ according to these amounts. Then through each edge before
 level $L_R(o)$ the amount that flows is at most the initial
 $k^{-R}$, while the capacity of such an edge is at least
 $\lambda^{-R}$, which is bigger because of (\ref{capacity}). So this
 is an admissible flow from $o$ to $L_R(o)$.

 The value of the flow at each vertex in $L_R(o)$ is at most
 $k^{-2R}$. For each such vertex $x$, we have $|L_R(x)|\geq
 (k-1)k^{R-1}$. Divide the amount at $x$ equally among these vertices, do
 the same for all $x\in L_R(o)$, and continue the flow from $L_R(o)$
 to $L_{2R}(o)$ according to this. Through each edge between $L_R(o)$
 and $L_{2R}(o)$ the amount that should flow is at  most $k^{-2R}$,
 while the capacity of an edge is at least  $\lambda^{-2R}$. Thus we
 have an admissible flow constructed already  from $o$ to $L_{2R}(o)$.

 Continuing in this manner: between the levels $L_{(n-1)R}(o)$ and
 $L_{nR}(o)$ the amount that should flow through an edge is at most
 $(k-1)^{-(n-2)}k^{-(nR-n+2)}$, and the capacity of such an edge is at
 least $\lambda^{-nR}$. This second is bigger because of
 (\ref{capacity}), thus we have constructed an admissible flow from
 $o$ to infinity with the positive value $k^{-R}$.\QED \medskip

 The converse is clearly false, as shown for example by a
 $(d+1)$-regular tree with an additional vertex placed on each of the
 edges. This tree $T_d'$ has branching number $\sqrt{d}$, while any
 vertex of the original $T_d$ together with its  $d+1$ neighbors form
 a $1$-fort of size $d+2$, so $p(T_d',2)=1$.  Furthermore, one might
 ask how sharp the above result $\br(T)\geq k-\frac{2k\log k}{\log N}$ is
 --- see Section \ref{open}.  \medskip

\noindent{\it Proof of Theorem \ref{pruning}.} Fix a vertex $x$ as the
root of $T$. It is enough to prove the theorem for $k=1$, and thus
find a small 1-fort $F_1$ in $T$ containing $x$, because then we can
inductively find a 1-fort $F_2$ inside $F_1$ with $\br(F_2)\leq
\br(F_1)-1 \leq \br(T)-2$, which will also be a 2-fort in $T$, and so
on.

By the Max-Flow-Min-Cut theorem, the branching number is characterized
by   \be\label{cut} \br(T)=\sup\left\{\lambda:\inf_{\Pi}
\sum_{e\in\Pi}\lambda^{-|e|}>0\right\}, \ee where the $\inf$ is over
cutsets $\Pi$ of edges separating $x$ from $\infty$. The expression
$\mu_\lambda(\Pi):=\sum_{e\in\Pi}\lambda^{-|e|}$ will be called the
{\bf $\lambda$-content} of the cutset (or of an arbitrary set of
edges).
\smallskip

Fix some $\beta>1$, and take an arbitrary finite tree $\T$ with
root $r$. By its boundary $\partial\T$ we mean the set of edges
with a leaf as an endpoint. If $\T=\{r\}$, then let
$\mu_{\beta}(\partial\T)=1$. Otherwise, denote the children of $r$
by $r_1,\dots,r_\ell$. Deleting the edge $(r,r_i)$ from $\T$
results in two connected components; the subtree that contains
$r_i$ will be denoted by $\T_i$, and $\T_i$ together with $r$ by
$\hat\T_i$. We have the disjoint union
$\cup_{i=1}^\ell\partial\hat\T_i=
\partial\T$, hence $m_1+\dots+m_\ell=\mu_\beta(\partial\T)$, where
$m_i:=\mu_\beta(\partial\hat\T_i)$.  We may assume $m_1\leq
m_2\leq\dots\leq m_\ell$. Now let us delete from $\T$ the entire
``$\beta$-largest'' subtree $\T_\ell$. Then look  at the subtrees
$\T_1,\dots,\T_{\ell-1}$, and repeat the whole procedure with each
$\T_i$ instead of $\T$, with root $r_i$, deleting the
``$\beta$-largest'' subtree from each $\T_i$. Repeat this
procedure over and over again until reaching the boundary of $\T$
in all subtrees. The remaining subtree $\F$ is clearly a 1-fort
inside $\T$. We claim that
\be\label{pruned}
\mu_{\beta-1}(\partial\F\cap\partial\T)\leq \mu_\beta(\partial\T)^\alpha,
\ee
where $\alpha=\beta/(\beta-1)$. Equality holds only for finite
$\beta$-ary trees $\T$, for integer $\beta$.

Before proving this claim, we show how it implies the existence of an
infinite 1-fort $F$ inside $T$, rooted at $x$, with
$\br(F)\leq\br(T)-1$.
\smallskip

Take a strictly decreasing sequence of positive numbers
$\{\beta_n\}$ converging to $\br(T)\geq 1$. Let
$\alpha_n=\beta_n/(\beta_n-1)>1$.
We can suppose that $T$ has no leaves. We have $\beta_1>\br(T)$,
so by the characterization (\ref{cut}), for any $\epsilon_1>0$
there exists a cutset $\Pi^1$ separating $x$ from $\infty$ with
$\mu_{\beta_1}(\Pi^1)<\epsilon_1$. If the finite subtree between
$x$ and $\Pi^1$ is called $\T^1$, then our above procedure finds a
1-fort $\F^1$ of $\T^1$, with
$\mu_{\beta_1-1}(\partial\F^1\cap\Pi^1)<\epsilon_1^{\alpha_1}$.
This upper bound is less than $1/2$ if we choose $\epsilon_1=1/2$.
Now denote the lower endvertices of the edges in
$\partial\F^1\cap\Pi^1$ by $x_1,\dots,x_\ell$.  The infinite
subtree of $T$ starting at $x_i$, called $T_i$, has branching
number less than $\beta_2$. Hence, for each $i$ and any
$\epsilon_2>0$, we can take a cutset $\Pi^2_i$ separating $x_i$
from $\infty$ with $\mu^*_{\beta_2}(\Pi^2_i)<\epsilon_2$, where
$\mu^*_{\beta_2}$ denotes $\beta_2$-content with distances $|e|$
measured from the new root $x_i$. That is,
$\mu_{\beta_2}(\Pi^2_i)<\epsilon_2\beta_2^{-|x_i|}$. If the finite
subtree between $x_i$ and $\Pi^2_i$ is called $\T^2_i$, then our
pruning procedure yields a 1-fort $\F^2_i$ inside each $\T^2_i$,
satisfying $\mu^*_{\beta_2-1}(\partial\F^2_i\cap\Pi^2_i)<
\epsilon_2^{\alpha_2}$. If we take the union
$\Phi^2:=(\partial\F^2_1\cap\Pi^2_1)\cup\dots\cup
(\partial\F^2_\ell\cap\Pi^2_\ell)$, then $\mu_{\beta_2-1}(\Phi^2)
= \sum_{i=1}^\ell \mu_{\beta_2-1}(\partial\F^2_i\cap\Pi^2_i)<
\sum_{i=1}^\ell \epsilon_2^{\alpha_2}(\beta_2-1)^{-|x_i|}=
\epsilon_2^{\alpha_2} \mu_{\beta_2-1}(\partial\F^1\cap\Pi^1)$.
Since $\mu_{\beta_2-1}(\partial\F^1\cap\Pi^1)$ is a finite number
independent of $\epsilon_2$, we can choose $\epsilon_2$ so small
that the last upper bound is less than $1/4$. Now we repeat
everything with the infinite subtrees of $T$ starting at the lower
endvertices $\{y_i, i=1,\dots,m\}$ of $\Phi^2$, using $\beta_3$
and some $\epsilon_3>0$. This gives a collection of cutsets
$\{\Pi^3_i\}$ and finite 1-forts $\{\F^3_i\}$. Take the union
$\Phi^3:=(\partial\F^3_1\cap\Pi^3_1)\cup\dots\cup
(\partial\F^3_m\cap\Pi^3_m)$, and choose $\epsilon_3$
sufficiently small so that $\mu_{\beta_3-1}(\Phi^3)<1/8$. Repeat this
ad infinitum, choosing $\epsilon_n$ such that
$\mu_{\beta_n-1}(\Phi^n)<2^{-n}$.

The union of all the finite 1-fort-pieces,
$F:=\F^1\cup\F^2_1\cup\dots\cup \F^2_\ell\cup\dots$, is an
infinite 1-fort of $T$, and each $\Phi^n$ is a cutset of $F$
separating $x$ from $\infty$. For any fixed $\beta>\br(T)$, if $n$
is large enough to have $\beta_n<\beta$, then
$\mu_{\beta-1}(\Phi^n)<\mu_{\beta_n-1}(\Phi^n)<2^{-n}$. Thus, by
definition (\ref{cut}), $\br(F)\leq \beta-1$. Since this holds for
all $\beta>\br(T)$, we have proved $\br(F)\leq \br(T)-1$.
\smallskip

We prove (\ref{pruned}) by induction on the depth of $\T$. If this
depth is 1, i.e.~each child $r_i$ of $r$ is a leaf, then $\F$ is just
obtained by deleting $r_\ell$, so we need to prove
$(\ell-1)/(\beta-1)\leq(\ell/\beta)^\alpha$. By taking derivatives
with respect to $\ell$ it is easy to check that the only value of
$\alpha$ for which this inequality holds for all real $\ell\geq 1$ is the
chosen $\alpha=\beta/(\beta-1)$. Equality holds only for $\ell=\beta$.

Suppose inductively that inside each subtree $\T_i$,
$i=1,\dots,\ell$, we have our 1-fort $\F_i$ with
$\mu^*_{\beta-1}(\partial\F_i\cap\partial\T_i)\leq m_i^\alpha$,
where $\mu^*_\lambda$ denotes  $\lambda$-content measured inside
$\T_i$ with root $r_i$, and $m_i=\mu^*_\beta(\partial\T_i)$. We
get $\F$ by joining the subtrees $\F_1,\dots,\F_{\ell-1}$ at $r$,
where $m_1\leq\dots\leq m_{\ell-1}\leq m_\ell$. Note that for
$\ell=1$ the claim is obvious. Now
$\mu_{\beta-1}(\partial\F\cap\partial\T)=
\left(\sum_{i=1}^{\ell-1}\mu^*_{\beta-1}(\partial\F_i\cap\partial\T_i)\right)/(\beta-1)
\leq\left(m_1^\alpha+\dots+m_{\ell-1}^\alpha\right)/(\beta-1)$,
while
$\mu_{\beta}(\partial\T)=\left(m_1+\dots+m_\ell\right)/\beta$.
Therefore, we would like to prove that
\be\label{ineq1}
\frac{m_1^\alpha+\dots+m_{\ell-1}^\alpha}{\beta-1} \leq
\left(\frac{m_1+\dots+m_\ell}{\beta}\right)^\alpha
\ee
for all possible values of the $m_i$'s. 
Let $y=(m_1+\dots+m_\ell)/m_\ell$. Because of $\alpha>1$ we have 
$m_i^\alpha\leq m_im_\ell^{\alpha-1}$. Adding these inequalities up, 
we get $(m_1^\alpha+\dots+m_{\ell-1}^\alpha)/(\beta-1)\leq
(y-1)m_\ell^\alpha/(\beta-1)$. Now recall that we proved 
$(y-1)/(\beta-1)\leq (y/\beta)^{\alpha}$ in the previous paragraph, 
hence our last upper bound is at most $(ym_\ell/\beta)^\alpha$. 
But this is just the RHS of (\ref{ineq1}), thus 
the proof of Theorem \ref{pruning} is complete.
\QED

\section{Regular graphs with anchored expansion}\label{anchored}

A simple generalization of the result \cite{sch92} for the Cayley
graph $\Z^2$ with standard generators is proved by \cite[Proposition
2.6]{grgr}: for any symmetric generating set of $\Z^2$, the
$2k$-regular Cayley graph $\Gamma_{2k}$ has  $p(\Gamma_{2k},k)=0$ and
$p(\Gamma_{2k},k+1)=1$.  As we have seen, the \crit probabilities for
regular trees all lie strictly between 0 and 1. Theorem
\ref{anchthm} suggests that this contrast between $\Z^d$ and the free
groups might have a geometric reason (see also the end of Section
\ref{open}). Indeed, the proof of the theorem will be based on the 
``perimeter method'', see in \cite{bape}.
\smallskip

\noindent{\it Proof of Theorem \ref{anchthm}.} Given an initial
configuration of occupied vertices, a set $S\subseteq V(G)$ is  called
{\bf internally spanned} if it becomes completely occupied even  in
the process restricted to $S$, i.e.~if we set all vertices in
$V(G)\setminus S$ to be vacant forever. First of all, we claim that if
complete occupation of $G$ occurs, then for any fixed vertex $o\in
V(G)$ there exists a strictly increasing sequence of finite connected 
internally spanned sets $o\in V_1\subset V_2\subset\dots\subset V(G)$.

If the vertex $o$ becomes occupied, then it does so in finite time, so
there exists a finite vertex set $V_1$ such that $o$ becomes occupied
even in the finite process restricted to $V_1$.  If we choose $V_1$ to
be minimal, then it is clearly a connected internally spanned set
containing $o$. Then let $V_1^*=V_1\cup\{v\}$ for some vertex $v$
neighboring $V_1$. Each vertex of $V_1^*$ becomes occupied in finite
time, so there is a minimal finite set $V_2$ such that all of $V_1^*$
becomes occupied even if the process is restricted to $V_2$. This
finite set $V_2$ is internally spanned, connected, and strictly larger
than $V_1$. Repeating this construction, we get the desired sequence
of random sets $o\in V_1\subset V_2\subset\dots$. Let us note that
with a bit more care one can achieve $\cup_{n\geq 1}V_n=V(G)$,
as well, but we will need only that there are arbitrarily large finite
connected internally spanned sets containing $o$.

Denote $v_n=|V_n|$ and $w_n=|\partial_eV_n|$, and take some
$0<h<\anch(G)$ such that $h+2k-d>0$ still holds.  The anchored
expansion property ensures that  $w_n/v_n\geq h$ for all sufficiently
large $n$.

Look at the $k$-neighbor process restricted to an internally
spanned $V_n$. If there are $x_n$ initially occupied vertices in
$V_n$, then the number of edges between these occupied vertices and
all the vacant vertices of $G$ (i.e.~the boundary of the occupied
part) is at most $dx_n$ initially. When a vacant vertex becomes
occupied, the boundary will have at most $d-k$ new edges, while at
least $k$ old edges disappear, so the boundary increases by at
most $d-2k$.  By the end of the complete occupation of $V_n$, we have
occupied $v_n-x_n$ initially vacant vertices, and have ended up
with a boundary $w_n$. Therefore,
$$dx_n+(v_n-x_n)(d-2k)\geq w_n>h v_n, \mbox{ so}$$
\be\label{proportion} x_n>v_n\frac{h-d+2k}{2k}=: cv_n.  \ee

Now take an i.i.d.~Bernoulli$(p)$ initial configuration on the whole infinite
graph, with $0<p<c$. Then, for any finite set $S\subset V(G)$,
\be\label{ld} \Pr_p(S\mbox{ contains at least }c|S|\mbox{ initially
occupied vertices})<e^{-I_p(c)|S|}, \ee  by the Large Deviation
Principle, see \cite[Theorem 2.1.14]{DZ}, where
$$I_p(c)=c\log\frac{c}{p} + (1-c)\log\frac{1-c}{1-p}\sim c\log\frac{1}{p}$$
when $c$ is fixed and $p\to 0$.

By a beautiful, by now well-known percolation argument from
\cite{kes}, in a $d$-regular graph there are at most $((d-1)e)^m$
possible connected sets $S\ni o$ (usually called ``lattice
animals'') of size $|S|=m$. Therefore, putting everything
together, for all large enough $M>0$,
\be\Pr_p(\mbox{complete occupation})&\leq&\Pr_p(\exists\mbox{  internally spanned }S\ni
o\mbox{ with }|S|>M)\nonumber\\
&\leq&\sum_{m=M}^{\infty}e^{-I_p(c)m+(\log(d-1)+1)m} \nonumber\\
&\to& 0,\mbox{ as }M\to\infty,\mbox{ if }I_p(c)>\log(d-1)+1.\nonumber
\ee

Thus $\Pr_p(\mbox{complete occupation})=0$ for $I_p(c)>\log (d-1)+1$,
which holds for all small enough $p>0$, in particular,
for $p<K(c)/(de-e)^{1/c}$, where $K(c)=c(1-c)^{(1-c)/c}$. Therefore,
$p(G_d,k)\geq K(c)/(de-e)^{1/c}>0$. 
\qed
\medskip

For the $(d+1)$-regular $T_d$, the above upper bound on the number of
lattice animals roughly coincides with the true asymptotics
$Cm^{-3/2}\left[d^d/(d-1)^{d-1}\right]^m$, see e.g.~\cite{pit98}.
However, for $d+1=2k$, $\anch(T_d)=d-1$, $c=\frac{d-1}{d+1}$,
the resulting estimate
$p(T_d,\lceil (d+1)/2\rceil)> \frac{1-o(1)}{(d-1)e}$ is very weak
compared to the true value $\sim 1/2$ coming from (\ref{middle}).

The sharpness of our theorem is shown by the free product $\Z^2*\Z$
with its natural 6-regular non-amenable Cayley-graph: from
$p(\Z^2,2)=0$ it follows immediately that $p(\Z^2*\Z,2)=0$. 
This also shows that the positivity result Proposition \ref{positive}
cannot be generalized to graphs with fast growth.

Theorem \ref{anchthm} can easily be used to give examples of non-trivial 
\crit probabilities for regular graphs that are not trees. For instance, 
the natural 4-regular Cayley graph of $\Z_3*\Z_3$ has no finite 1-forts,
so by Lemma \ref{nofinite} we get $0<p(\Z_3*\Z_3,2)<1$. For a more 
general result, see the end of Section \ref{open}.

\section{Concluding remarks and open problems}\label{open}

The Red Lemma \ref{red} and the Monotonicity Lemma \ref{monotone} 
give that having no finite $(k-1)$-forts implies $p(T,k)\leq
p(T_k,k)=1-1/k$. There are examples showing that, in general,
having no $(k-1)$-forts with $\br(F)<b-k+1$, where $b>k$ integer,
does not imply $p(T,k)\leq p(T_b,k)$. On the other hand, the
1-fort $F$ found by Theorem \ref{pruning} is the largest possible
(in terms of the $e^{\text{br}(T)-1}$-dimensional Hausdorff
measure of the boundary space $\partial F$, where the distance
between two infinite rays $\xi,\eta\in\partial F$ is
$e^{-|\xi\wedge\eta|}$, see \cite{lpbook}) when $T$ is a
$\br(T)$-regular tree. This suggests that regular trees might play
the role of extreme cases in the sense that $f_k(b)=p(T_b,k)$ for
all $b\in\N$ for the function in (\ref{function}), i.e.~they might
be the trees with a fixed branching number which are the easiest
to occupy. This would also nicely coincide with similar results
for random walks, see \cite{vir:speed} and \cite{vir:fast}.
However, as we will show below, for $b>k$ this is not the case,
even for Galton-Watson trees, for which Theorem \ref{pruning}
holds even with random pruning. So we are left with the following
open problem:
\smallskip

\noindent{\bf The easiest trees to occupy.} {\it Determine the
function $f_k(b)$. Is it strictly positive for all real $b\geq 1$?
Is it continuous apart from $b=k$?}
\smallskip

It is possible that $f_k(k)=1-\frac{1}{k}$. Also note that requiring a 
fixed bound on the degrees instead of the branching number already 
implies strict positivity, by Proposition \ref{positive}. 
\medskip

\noindent{\bf Galton-Watson trees.} One can study the same problems on
a Galton-Watson tree $T_\xi$ with offspring distribution $\xi$.  For
any $p$, the event  $\left\{\Pr_p(\mbox{complete occupation of
}T_\xi)>0\right\}$ is an inherited event, so it has probability 0 or
1, see \cite[Proposition 4.6]{lpbook}, which shows that $p(T_\xi,k)$
is a constant almost surely, given non-extinction. If $\Pr(\xi<k)>0$,
then infinitely many  finite $(k-1)$-forts of bounded size occur, so
$p(T_\xi,k)=1$. Otherwise, $T_\xi$ can be  built up from copies of
$T_k$, and we get $p(T_\xi,k)\leq p(T_k,k)=1-1/k$.
We also have $\Pr_p(\mbox{complete occupation of }T_\xi)=1$ a.s.,
given non-extinction, for $p>p(T_\xi,k)$. Just as above, this shows the
following  monotonicity property. If two offspring distributions $\xi$
and $\eta$ satisfy $\Pr(\xi<m)\geq\Pr(\eta<m)$ for all $m=1,2,\dots$,
i.e.~$\eta$ stochastically dominates $\xi$,
then there is natural coupling between the trees $T_\xi$ and $T_\eta$ such 
that $T_\eta$ is uniformly bigger than $T_\xi$ a.s., and so we get 
$p(T_\xi,k)\geq p(T_\eta,k)$.
\medskip

\noindent{\bf A GW tree beating a regular tree.} Consider the GW tree
$T_{\xi}$ with root $r$ and offspring distribution
$\Pr(\xi=2)=\Pr(\xi=4)=1/2$. Then $\br(T_\xi)=\E\xi=3$ a.s.~\cite{ly:perc90},
there are no finite $1$-forts in $T_\xi$, and $0<p(T_\xi,2)<1$ is an
almost sure constant. We claim that $p(T_\xi,2)<p(T_3,2)=1/9$.

Let $\RR(x,T_\xi)$ be the event $\{$the vertex $x$ of $T_\xi$ is in an
infinite vacant 1-fort$\}$, and set $q(T_\xi)=\Pr_p(\RR(r,T_\xi))$.
This is not an almost sure constant, so let us take expectation over
all GW trees: $q=\E(q(T_\xi))$. Now
$$q=\frac{1}{2}\E(q(T_\xi)\,|\, \xi_r=2)+ \frac{1}{2}\E(q(T_\xi) \,|\,
\xi_r=4).$$ Regarding the first term, $\Pr_p(\RR(r,T_\xi)\,|\,
\xi_r=2)=\Pr_p\bigl(r$ is initially vacant, and at least one of
$\RR(r_1,T_\xi')$ and $\RR(r_2,T_\xi'')$ does not fail$\bigr)$,
where $r_1$, $r_2$ are the two children of $r$, and $T_\xi'$,
$T_\xi''$ are the corresponding subtrees.
By the independence of initial configurations in
$T_\xi'$ and $T_\xi''$, this is equal to $(1-p)\bigl(
\Pr_p(\RR(r_1,T_\xi')) + \Pr_p(\RR(r_2,T_\xi'')) -
\Pr_p(\RR(r_1,T_\xi'))\Pr_p(\RR(r_2,T_\xi''))\bigr)$. Now, by the
recursive structure of $T_\xi$, and the independence of the subtrees
$T_\xi'$ and $T_\xi''$, taking the conditional expectation gives
$$\E(q(T_\xi)\,|\, \xi_r=2)=(1-p)(2q-q^2).$$ A similar argument for
the second term gives
$$\E(q(T_\xi) \,|\, \xi_r=4)=(1-p)(4q^3-3q^4).$$

Altogether, we have the equation $q=\frac{1}{2}(1-p)(2q-q^2+4q^3-3q^4)$, and
need to determine the infimum of $p$'s for which there is no solution
$q\in (0,1]$ --- that infimum will be $p(T_\xi,2)$. Setting
$f(q)=2-q+4q^2-3q^3$, an examination of $f'(q)$ gives that
$\max\{f(q):q\in [0,1]\}=f((4+\sqrt{7})/9)=2.2347\dots$. So there is
no solution $q>0$ if{f} $2/(1-p)>2.2347\dots$, which gives
$p(T_\xi,2)={0.10504\dots}< 1/9$.\QED
\medskip

\noindent{\bf Small $k$-forts.}  If we define $\Gamma_k(N)$ as the set
of trees without $k$-forts of size at most $N$, and
$\gamma_k(N)=\inf\left\{\br(T) : T\in\Gamma_k(N)\right\}$, then we
know only \be\label{extremal} 2 - 2^{-(1+o(1))N/2} \geq \gamma_1(N)
\geq  2-\frac{c}{\log N}.  \ee The upper bound is achieved by trees
analogous to the tree on the left in Figure \ref{1-forts}. (The
vertices of degree two are distributed following a greedy strategy: let 
the root be the first one of them, and then, at the generic step, put
them on the highest level possible, in the highest possible number
at that level, subject to not forming a 1-fort of size at most $N$.)
Actually, this gives asymptotically the smallest branching number that
a tree in $\Gamma_1(N)$ with maximal degree $3$ can have. The lower
bound comes from the proof of Theorem \ref{maintree}. 
\medskip 

\noindent{\bf Amenable and non-amenable groups.}
As we already discussed in Section \ref{anchored}, for the free
Abelian groups $\Z^d$ the \crit probabilities are
almost completely determined by \cite{sch92} and \cite{grgr}.
The simplest non-Abelian group, the {\bf Heisenberg group}, can be
considered with natural generator sets of 2 or 3 elements \cite{dlH},
and it seems reasonable to conjecture that the corresponding 4- or
6-regular Cayley graphs $H_4$ and $H_6$ have $p(H_{2k},k)=0$.
One can easily find finite $k$-forts to prove $p(H_{2k},k+1)=1$.

The most famous amenable groups with exponential growth are the
{\bf lamplighter groups} $\Z_r \wr \Z^d$.  With a natural
generating set, the Cayley graph of $\Z_r \wr \Z$ is the {\bf
Diestel-Leader graph} $DL(r,r)$, where $DL(r,s)$ is the
``horocyclic product'' of two regular trees $T_r$ and $T_s$, see
\cite{DLlamp}. These transitive graphs with degree $r+s$ are
amenable if{f} $r=s$, and it is conjectured that for $r\not=s$
they are not quasi-isometric to any Cayley graph. It is not
difficult to see that k-neighbor bootstrap percolation on
$DL(r,s)$, where $r\leq s$, has critical probability $1$ for
$k>s$, while strictly between $0$ and $1$ for $r+1\leq k\leq s$,
if such $k$ exists. However, it is unclear if $p(DL(r,s),r)=0$
holds or not. A positive answer, together with our proof of
Theorem \ref{anchthm}, would have the interesting consequence
that, as $p$ gets closer and closer to 0, complete occupation of
$DL(r,r)$ by the $r$-neighbor rule will happen more and more
through F\o lner sets, rather than through the exponentially
growing balls.

If a finitely generated non-amenable group $G$ contains a free
subgroup on two elements, then, as David Revelle pointed out, there
exists a generating set of $k$ elements with $0<p(G_k,\lceil
k/2\rceil)<1$. The reason is that if $G$ is originally defined by a
symmetric generating set of $t$ elements, then taking $2t$ free
symmetric generators inside the free subgroup, we arrive at a
$k=t+2t$-regular graph, in which each vertex is contained in a
$2t$-regular subtree. So our results give
$0<p(G_k,\lceil 3t/2 \rceil)\leq p(T_{2t-1},\lceil 3t/2 \rceil)<1$.
\smallskip

\noindent{\bf An open question} inspired by the above results:
{\it is a group amenable if and only if for any finite generating set,
the resulting $k$-regular Cayley graph has $p(G_k,\ell)\in\{0,1\}$ for
any $\ell$-neighbor rule?}
\medskip

\noindent{\bf Acknowledgments.} We are grateful to Dayue Chen,
Manjunath Krishnapur, Fabio Martinelli, Robin Pemantle, David
Revelle, \'Ad\'am Tim\'ar, B\'alint Vir\'ag and the referee 
for helpful discussions and comments.

\end{document}